\documentclass{article}
\usepackage{amsmath}
\usepackage{amssymb}
\usepackage[dvipdfmx]{graphicx}
\usepackage{mathrsfs} 
\usepackage{tabu}
\newtheorem{theorem}{Theorem}
\newtheorem{lemma}{Lemma}

\newtheorem{remark}{Remark}
\newcommand{\mbf}[1]{\mbox{\boldmath{${#1}$}\unboldmath}}
\begin{document}
\title{Smoothed nonparametric two-sample tests}

\author{
Taku MORIYAMA \and Yoshihiko MAESONO
}
\date{}

\maketitle

\begin{abstract}
We propose new smoothed median and the Wilcoxon's rank sum test. As is pointed out by Maesono et al.\cite{maesono2016smoothed}, some nonparametric discrete tests have a problem with their significance probability. Because of this problem, the selection of the median and the Wilcoxon's test can be biased too, however, we show new smoothed tests are free from the problem. Significance probabilities and local asymptotic powers of the new tests are studied, and we show that they inherit good properties of the discrete tests. 
\end{abstract}

Let $X_1,X_2,\cdots,X_m$ be independently and identically distributed random variables ({\it i.i.d.}) from distribution function $F(x)$ and $Y_1,Y_2,\cdots,Y_n$ be {\it i.i.d.} from $F(x-\theta)$ where $\theta$ is unknown location parameter. We assume that $m+n =N$ and $\lambda_N := m/N \to \lambda$ s.t. $0< \lambda <1$. We consider `2-sample problem' whose null hypothesis $H_0 : \theta=0$ and alternative $H_1 : \theta >0$. There are many nonparametric tests based on linear order statistic (see H\'{a}jek et al.\cite{sidak1999theory}). The median and the Wilcoxon's rank sum test are widely used and investigated well among them. Moreover, the median test is known for its low cost especially in survival analysis because it only needs the number of the lower values than the combined sample median. However, Freidlin \& Gastwirth\cite{freidlin2000should} have reported that the median test is numerically inferior even in double exponential case in spite of its theoretical most powerfulness. We first investigate the median's power again and find theoretical and numerical superiority when the distribution has heavy tail.

As is shown by Maesono et al.\cite{maesono2016smoothed}, the sign test and the Wilcoxon's signed rank test have a problem of their significance probabilities, and it comes from their discreteness of $p$-values. Because of this, the discrete tests can let an distorted statistical decision, and we show that the median and the Wilcoxon's rank sum test are also the same. Using a smooth function in a similar way as Maesono et al.\cite{maesono2016smoothed}, we propose the smoothed test statistics of them, and the smoothed median test also needs only the lower values than the combined median. It is proved that their Pitman's $A.R.E.$ are respectively same, and they are asymptotically nonparametric.

\section{The median and the Wilcoxon's rank sum test}

\subsection{The median and the Wilcoxon's rank sum test}

In this section, we introduce the median and the Wilcoxon's rank sum test statistic and their properties. Some numerical experiments on their power are studied and the problem about their $p$-values is brought up here.

Let us define that $\psi(x)=1 ~(x\geq0), ~=0 ~(x<0)$ and $Z_{1}=X_{1}, \cdots, Z_{m}=X_{m}, \cdots, Z_{m+1}=Y_{1}, \cdots, Z_{N}=Y_{n}$. There are various forms of the median test statistic, and here we define
$$M=M({\mbf X},{\mbf Y})=\sum_{j=1}^n\psi(Y_j -Z)$$
where ${\mbf X}=(X_1,X_2,\cdots,X_m)^T$ and $Z$ denotes the sample median of $\{Z_1,Z_2,\cdots,Z_N\}$. We put $Z=Z_{((N+1)/2)}$ if $N$ is odd and $Z=(Z_{(N/2)} +Z_{((N/2)+1)})/2$ else where $Z_{(1)} < \cdots < Z_{(N)}$ are the order statistics. The Wilcoxon's rank sum test statistic is given by
$$W_{2}=W_{2}({\mbf X},{\mbf Y})=\sum_{1\leq i\leq m} \sum_{1\leq j\leq n}\psi(Y_j-X_i).$$
For observed values ${\mbf x}=(x_1,x_2,\cdots,x_m)^T$ and ${\mbf y}=(y_1,y_2,\cdots,y_n)^T$, we put $m=m({\mbf x},{\mbf y})$ and $w_2=w_2({\mbf x},{\mbf y})$ which are the realized values of $M$ and $W_2$.  If the $p$-value
$$P_0(M\geq m)\hspace{5mm}{\rm or} \hspace{5mm}P_0(W_2\geq w_2)$$
is small enough, we reject the null hypothesis $H_0$.  Here $P_0(\cdot)$ denotes a probability under the null hypothesis $H_0$.

As we can see easily, the median test only needs the number of the lower values than the combined median, and we can finish the observation when the median is obtained. If the data size is large or the tail of the distribution is heavy, it is possible to save much cost and time. However, not many statisticians take the great merit into account. We show the smoothed median test inherits the superiority in Section 3.1.

\subsection{Local power}

The median and the Wilcoxon's rank sum test are one of linear rank tests which are defined as follows
$$S=S({\mbf X},{\mbf Y}) =\sum_{j=1}^n a_{m,n}(R_j)$$
where $R_j (j=1,\cdots,n)$ denotes rank of the observation $Y_j$ in the combined sample $Z_{1}, \cdots, Z_{N}$. The Wilcoxon's test is given by 
$$a_{m,n}(u) = m+ \frac{1}{2}(n+1) -u,$$
so the test statistic is
$$W_{2}=W_{2}({\mbf X},{\mbf Y})= mn+ \frac{1}{2}n(n+1) -\sum_{j=1}^n R_j.$$
The last term
$$U=\sum_{j=1}^n R_j$$
is led as locally most powerful rank test in logistic case, and $W_2$ is equivalent to the test statistic $U$ (H\'{a}jek et al. \cite{sidak1999theory}). The median test statistic is obtained by
$$a_{m,n}(u) = \psi(u- \frac{1}{2}(m+n+1)),$$
and so
$$M=M({\mbf X},{\mbf Y})=\sum_{j=1}^n \psi(R_j - \frac{1}{2}(m+n+1)).$$
$M$ is asymptotically equivalent to the locally most powerful rank test statistic when the underlying distribution is the double exponential.

Table 1 shows Pitman's asymptotic relative efficiencies ($A.R.E.$) of the combinations of the two-sample t-test $T_2$, the median and the Wilcoxon's test. Pitman's $A.R.E.$s are given by a ratio of the values related to asymptotic local power of the two tests. We confirm that $T_2$ is the most powerful in normal case, $M$ is in double exponential case and so on. Note that $ARE(M|W_2) = 1.33$ in the double exponential case.

\begin{table}[h]
\caption{Pitman's $A.R.E.$ of 2-sample tests \label{tab:table1}}
\begin{center}
\begin{tabular}{c|ccc}\hline
population distribution&Nor.&Logis.&D.exp.\\\hline
$ARE(M|T_2)$&0.637&0.822&2\\
$ARE(W_2|T_2)$&0.955&1.10&1.5\\
$ARE(M|W_2)$&0.667&0.75&1.33\\\hline
\end{tabular}
\end{center}
\end{table}

Nevertheless, the numerical weakness of the median test's power was reported. Freidlin \& Gastwirth \cite{freidlin2000should} shows the empirical local power of some two-sample linear rank tests, and in their tables, the median test is inferior to the Wilcoxon's test even in the double exponential case.

We investigate their powers both theoretically and numerically in heavy tailed case. Table 2 shows Pitman's $ARE(M|W_2)$ of $T(2)$ ($T$ distribution with $2$ degrees of freedom), $T(1)$ (Cauchy distribution) and $T(1/2)$. Distribution's tail gets heavy in accordance with decrease in the degree, and we can prove the monotonic increase of $ARE(M|W_2)$. 

\begin{table}[h]
\caption{Pitman's $A.R.E.$ of 2-sample tests in $T$ distribution \label{tab:table2}}
\begin{center}
\begin{tabular}{c|ccc}\hline
population distribution&$T(2)$&$T(1)$&$T(1/2)$\\\hline
$ARE(M|W_2)$&0.961&1.33&2.29\\\hline
\end{tabular}
\end{center}
\end{table}

Table 3 and 4 are the numerical results of the ratio of their empirical local power in the $t$ distributions. We find that the power of the median test is rather stronger especially in $T(1/2)$. If to obtain complete data needs much cost and time in such case, the median test has great superiority.

\begin{table}[h]
\caption{Ratio of empirical local power in $T(1)$ distribution \label{tab:table2}}
\begin{center}
\begin{tabular}{c|cccccc}\hline
sample size ($m,n$) &(10,10)&(20,20)&(30,30)&(50,50)&(10,30)&(30,10)\\\hline
$\theta =1$&1.03&1.16&1.17&1.09&1.06&1.06\\
$\theta =0.5$&1.01&1.15&1.17&1.16&1.07&1.06\\
$\theta =0.1$&0.994&1.03&1.05&1.05&1.04&0.993\\\hline
\end{tabular}
\end{center}
\end{table}

\begin{table}[h]
\caption{Ratio of empirical local power in $T(1/2)$ distribution \label{tab:table2}}
\begin{center}
\begin{tabular}{c|cccccc}\hline
sample size ($m,n$) &(10,10)&(20,20)&(30,30)&(50,50)&(10,30)&(30,10)\\\hline
$\theta =1$&1.26&1.52&1.54&1.54&1.34&1.34\\
$\theta =0.5$&1.17&1.37&1.42&1.57&1.27&1.26\\
$\theta =0.1$&1.06&1.08&1.10&1.14&1.07&1.06\\\hline
\end{tabular}
\end{center}
\end{table}

\subsection{Significance probability}

We will show that the median test and the Wilcoxon's rank sum test have a problem with their significance probabilities. Maesono et al.\cite{maesono2016smoothed} reports that the sign and the Wilcoxon's signed rank test can make `distorted' statistical results, and this problem comes from their $p$-values' discreteness. The median and the Wilcoxon's test are also discrete, and we study their significance probabilities. 

Table 5 shows the ratio of frequency of exact $p$-value of $W_2$ smaller than that of $M$ in the following tale area $\Omega_{\alpha}$
$$\Omega_{\alpha}=\left\{{\mbf x}\in{\mbf R}^n ~\left|~ \frac{m({\mbf x})-E_0(M)}{\sqrt{V_0(M)}} > v_{1-\alpha},\right.\hspace{4mm}{\rm or}\hspace {5mm}\frac{w_2({\mbf x})-E_0(W_2)}{\sqrt{V_0(W_2)}} > v_{1-\alpha}\right\}$$ 
where $v_{1-\alpha}$ is a $(1-\alpha)th$ quantile of the standard normal distribution $N(0,1)$, and $E_0(\cdot)$ and $V_0(\cdot)$ stand for an expectation and a variance under $H_0$, respectively.  We count samples that an exact $p$-value of the test is smaller than the other in $\Omega_{\alpha}$, and calculate the ratio of the frequency.

\begin{table}[h]
\caption{Comparison of significance probabilities \label{tab:table5}}
\begin{center}
\begin{tabular}{c|cccccc}\hline
sample size ($m,n$)&(10,10)&(20,20)&(30,30)&(10,20)&(20,10)&($U_m^*,U_n^*$)\footnotemark[1]\\\hline
$z_{0.9}$&2.21&2.79&1.68&7.13&7.07&3.86\\
$z_{0.95}$&6.37&5.66&2.17&3.44&3.41&4.21\\
$z_{0.975}$&3.05&2.80&3.72&16.6&14.7&4.10\\
$z_{0.99}$&33.7&7.54&1.56&6.05&5.29&4.11\\\hline
\end{tabular}
\end{center}
\end{table}
\footnotetext[1]{$U_m^*$ and $U_n^*$ are random numbers from the discrete uniform distribution $U^*(5,40)$}

Because the values in Table 5 are larger than 1, we find that $W_2$ tends to have smaller $p$-value than $M$, and they can let us to use $W_2$ if we wants the small $p$-value. This comes from that the possible $p$-values of $M$ are more sparse than $W_2$ like the sign and Wilcoxon's signed rank test as is discussed by Maesono et al.\cite{maesono2016smoothed}.

In order to conquer the problem, we propose the smoothed median test $\widetilde{M}$ and the  Wilcoxon's test $\widetilde{W_2}$. The discrete tests are distribution-free, but the smoothed tests are not. However, we can confirm that they are asymptotically distribution-free and their Pitman's efficiencies of them are the same.

\section{Smoothed median test}

\subsection{Smoothed median test}

Hereafter, we assume that $N$ is odd for the brevity, and we consider to make $M$ smooth appropriately. A possible way is to define them as kernel type statistics, and we introduce a kernel distribution estimator first. The empirical distribution function of a population distribution function $F$ is given by
$$F_n(x) = \frac{1}{n} \sum_{i=1}^n\psi(x -X_i).$$
Let $k(u)$ be a kernel function which satisfies
$$\int_{-\infty}^{\infty}k(t)dt=1,$$
and $K(t)$ is an integral of $k(t)$ s.t.
$$K(t)=\int_{-\infty}^{t}k(u)du.$$
In this paper, we assume that the kernel $k$ is a symmetric function around the origin. The kernel distribution estimator of $F$ is given by
$$\widehat{F}(x)=\frac{1}{n}\sum_{i=1}^n K\Bigl(\frac{x-X_i}{h}\Bigr)$$
where $h$ is a bandwidth which satisfies $h\to 0$ ($n \to\infty$).

The median test statistic is given by
$$M = \frac{n-m+1}{2} -\sum_{i=1}^m\psi^*(Z -X_i)$$
where $\psi^*(x)=1 ~(x>0), ~=0 ~(x\leq0)$, and we put $M^{\dagger}$ as the second term. Applying the kernel smoothing to $M^{\dagger}$, we define the following smoothed test statistic
$$\widetilde{M}=\sum_{i=1}^m K^*\left(\frac{Z -X_i}{h}\right)$$
where $K^*$ is an integral of the kernel $k^*(t)$ which satisfies
\begin{eqnarray*}
\int_{0}^{\infty}k^*(t)dt=1 ~~~{\rm and}~~~ k^*(t)=0 ~~{\rm for}~~ t \le 0
\end{eqnarray*}
and $h$ is a bandwidth. In addition, we assume that $A_{1,1}^*=0$ where
$$A_{i,j}^* = \int_{-\infty}^{\infty} t^i {k^*}^j(t) dt.$$
As we see easily from the definition, the smoothed median test does not need the values of the data larger than the combined median.

\begin{remark}
The above condition of $k^*$ is not tough, and we can easily  construct the following simple polynomial type one
\begin{equation*}
k^*(t)=\left(-6t+4\right) I(0<u<1).
\end{equation*}
We can easily construct a kernel function which satisfies $A_{1,1}^*=A_{2,1}^*=0$ too. It is not a problem that $k^*$ may take negative value, because $\widetilde{M}$ is test statistic.
\end{remark}

\subsection{Asymptotic properties}

The median test statistic $M^{\dagger}$ exactly follows the hypergeometric distribution $HG(N,m,(N-1)/2)$ under $H_0$, and so we easily find that
$$E_{0}[M^{\dagger}] = \frac{m}{2}\left( 1-\frac{1}{m+n}\right),~~~ V_{0}[M^{\dagger}] =\frac{mn}{4(m+n)}.$$
Using the following joint distribution of $Z$ and $U (\in \{1,\cdots,r\})$ which is the number of $\{X_1, \cdots, X_m\}$ less than $Z$
\begin{eqnarray}
&&h(u,z) \nonumber\\
&=& m \binom{m-1}{u} \binom{n}{r-u}F^u(z)[1-F(z)]^{m-u-1}F^{r-u}(z- \theta)[1-F(z-\theta)]^{n-r+u} f(z) \nonumber\\
&&+ n \binom{m}{u} \binom{n-1}{r-u}F^u(z)[1-F(z)]^{m-u}F^{r-u}(z- \theta)[1-F(z-\theta)]^{n-r+u-1} f(z-\theta) \nonumber
\end{eqnarray}
under $H_1$, Mood \cite{mood1954asymptotic} obtains the following asymptotic expectation
$$E_{\theta}[M^{\dagger}]= mF(z_{\theta,N}) + o(N)$$
where $z_{\theta,N}$ stands for the median of the distribution $G_{\theta,N}$ defined as
$$G_{\theta,N}(x)= \lambda_N F(x) +(1-\lambda_N) F(x-\theta).$$
Then we find the following Pitman efficiency is given by
$$e_P[M^{\dagger}] = \lim_{N}\left[(N V_{0}[M^{\dagger}])^{-1/2} \frac{\partial}{\partial \theta} E_{\theta}[M^{\dagger}] \biggm|_{\theta=0} \right]= 2\sqrt{\lambda(1-\lambda)} f(z_0)$$
where $z_0$ is the median of the population distribution of $F$.

Now, we clarify the difference of the asymptotics between the odd and even. If $N$ is odd, the exact distribution of $Z$ is given by
$$F_{Z}^{O}(z) = \frac{1}{\beta((N+1)/2,(N+1)/2)} \int_{-\infty}^z [F(x)]^{(N-1)/2} [1-F(x)]^{(N-1)/2} f(x) dx$$
where $\beta(\cdot,\cdot)$ is the beta function. For an even number, the exact distribution is (Desu \& Rodine \cite{mahamunulu1969estimation})
$$F_{Z}^{E}(z) = \frac{2}{\beta(N/2,N/2)} \int_{-\infty}^z [F(x)]^{(N/2)-1} \Bigm\{ [1-F(x)]^{(N/2)} - [1-F(2z -x)]^{(N/2)}\Bigm \}f(x) dx.$$
From the above, we have the following lemma.
\begin{lemma}
The difference of the distribution of $Z$ between $N=\widetilde{N}$ and $N=\widetilde{N}+1$ is
$$| F_{Z}^{O}(z) -F_{Z}^{E}(z) | = o(\widetilde{N}^{-2}) ~~~ \text{(negligible)}.$$
\end{lemma}
\noindent{\it Proof.} The direct computation gives the result. See Appendices.

Hereafter, we denote the $l$-th derivative of $f$ by $f^{(l)}$. About the Pitman's efficiency, we obtain the following result.
\begin{theorem}
Let us assume that $f^{(1)}$ exists and is continuous at a neighborhood of both $z_0$ and $z_{\theta,N}$, and both $f(z_0)$ and $f(z_{\theta,N})$ are positive. In addition, we assume that $h=o(N^{-1/2})$ or that $A_{1,1}^*=0$ and $h=o(N^{-1/4})$. Then, we have
$$e_P[\widetilde{M}] = e_P[M^{\dagger}].$$
\end{theorem}
\noindent{\it Proof.} See Appendices.

Since the main term of $\widetilde{M}$ is a two-sample $U$ statistic, it is easy to prove the following asymptotic normality.
\begin{theorem}  Let us assume that $f^{(1)}$ exists and is continuous at a neighborhood of $z_0$, and $f(z_0) >0$. If $h=o(N^{-1/2})$, or $A_{1,1}^*=0$ and $h=o(N^{-1/4})$ holds, we have
$$\sup_{-\infty < x < \infty} \left| P_0 \left[ V_{1}^{-1/2} (\widetilde{M}-E_1) < x\right] -\Phi(x) \right| = o(1) $$
where
$$E_1=\frac{m}{2}\left(1- \frac{1}{m+n}\right),~~~ V_1=\frac{mn}{4(m+n)},$$
$\Phi$ is the standard normal distribution function and $\epsilon$ is any positive number.
\end{theorem}
\noindent{\it Proof.} The evaluation of the difference between $\widetilde{M}$ and $M^{\dagger}$ is the main, and the detail is in  Appendices.

Note that the main terms of the asymptotic expectation and variance of $\widetilde{M}$ do not depend on $F$ under $H_0$.

Next, we study the local power of $\widetilde{M}$, and obtain the following result.
\begin{theorem}  Under the same assumptions of Theorem 2.3 and that $\sqrt{N} h \to \infty$, if $A_{1,1}^*=0$ and $A_{1,1,1}^*$ is positive, we have
\begin{eqnarray*}
\lim_{N \to \infty} \frac{1}{h} (LP_{\frac{\xi}{\sqrt{N}}, \alpha}[\widetilde{M}] - LP_{\frac{\xi}{\sqrt{N}}, \alpha}[M^{\dagger}]) > 0.
\end{eqnarray*}
where
$$LP_{\frac{\xi}{\sqrt{N}}, \alpha}[\widetilde{M}] = P_{\frac{\xi}{\sqrt{N}}}\left[ {V_{0}[\widetilde{M}]}^{-1/2}(\widetilde{M} -E_0[\widetilde{M}]) > v_{1-\alpha}\right]$$
and $LP_{\frac{\xi}{\sqrt{N}}, \alpha}[M^{\dagger}]$ is also the same.
\end{theorem}
\noindent{\it Proof.} See Appendices.

\begin{remark}
We can construct the following polynomial type kernel
\begin{equation*}
k^*(t)=\left[\left( \frac{\pm 3\sqrt{4353} -3}{17}\right) u^2+ \left( \frac{\mp 3\sqrt{4353} -99}{17} \right) u + \left(\frac{\pm \sqrt{4353} +135}{34} \right)\right] I(0<u<1).
\end{equation*}
which satisfies $A_{1,1}^*=0$ and $A_{1,1,1}^*=1 (>0)$. Similarly, we have the following exponential type
\begin{eqnarray*}
k^*(t)&=&\biggm[e^{-t} + \left(\frac{613 -2\sqrt{207586}}{58} \right)*(2e^{-2t})+\left(\frac{3\sqrt{207586} -1137}{58} \right)*(3e^{-3t}) \\
&& ~~~ +\left(\frac{524 -\sqrt{207586}}{58} \right)*(4e^{-4t})\biggm]
\end{eqnarray*}
which satisfies $A_{1,1}^*=0$ and $A_{1,1,1}^*=1 (>0)$. We want to use a kernel whose value of $A_{1,1,1}^*$ is larger, however, in practice, the estimated value of $V_{1}^{\star}$ which we need can be negative.
\end{remark}

\section{Smoothed Wilcoxon's rank sum test}

\subsection{Smoothed Wilcoxon's rank sum test}

Here, we give the smoothed test $\widetilde{W}_2$ in the same manner. In the same way, we can define the smoothed test statistics of $W_2$
$$\widetilde{W}_2=\sum_{i=1}^{m} \sum_{j=1}^n K\left(\frac{Y_j-X_i}{h}\right).$$

The following moments of the Wilcoxon's test statistic $W_2$ are easy to obtain
$$E_{\theta}[W_2] = mn \int_{-\infty}^{\infty} f(y) F(y+\theta) dy,~~~ V_{0}[W_2] =\frac{mn(m+n)}{12},$$
and then we find the following Pitman efficiency
$$e_P[W_2] = \sqrt{ 12\lambda(1-\lambda)} \int f^2(x) dx.$$

Using variable changes and the Taylor expansion, we obtain the following asymptotic expectation of $\widetilde{W}_2$
\begin{eqnarray*}
E_{\theta}[\widetilde{W}_2] &=& mn \int_{-\infty}^{\infty} \int_{-\infty}^{\infty} K\left(\frac{x -y}{h}\right) f(x-\theta) f(y)dx dy  \\
&=& mn \int_{-\infty}^{\infty} \int_{-\infty}^{\infty} k\left(v\right) f(y +hv -\theta) F(y)dv dy\\
&=& mn \left[\int_{-\infty}^{\infty} f(y) F(y+\theta) dy +O(h^2) \right].
\end{eqnarray*}

Hereafter, we assume that $h = o(N^{-1/2})$ or that $A_{1,1}=\cdots=A_{3,1}=0$ and $h = o(N^{-1/4})$ in order to ignore the residual term. It is easy to obtain such kernels using Jones \& Signorini \cite{jones1997comparison}.

Further, it is easy to see that $V_{0}[\widetilde{W}_2]= mn(m+n)/12 + o(N^3)$. Combining the above results, we can obtain the following result.

\begin{theorem}
Let us assume that $f^{(1)}$ exists and is continuous at a neighborhood of both $0$ and $\theta$. Then, we have
$$e_P[\widetilde{W}_2] = e_P[W_2].$$
\end{theorem}
\noindent{\it Proof.} See Appendices.

Using the asymptotic results of two-sample $U$-statistics, we can obtain the following error bound of the normal approximation.
\begin{theorem}  Let us assume that $f^{(1)}$ exists and is continuous around a neighborhood of $0$ and $h = o(N^{-1/2})$, or that such $f^{(3)}$ exists, $A_{2,1}=0$ and $h = o(N^{-1/4})$. Then, we have
$$\sup_{-\infty < x < \infty} \left|P_0 \left[ V_2^{-1/2} (\widetilde{W}_2 -E_2) < x\right] -\Phi(x) \right| = o(N^{-1/2}) $$
where
$$E_2=\frac{mn}{2},~~~ V_2=\frac{mn(m+n)}{12}.$$
\end{theorem}
\noindent{\it Proof.} We need the theory of edgeworth expansion of the  two-sample $U$-statistics. See Appendices.

In the same manner as the local power of $\widetilde{M}$, we can obtain the power of $\widetilde{W}_2$, and find the following result.

\begin{theorem}  Under the assumptions of Theorem 3.1, we have
\begin{eqnarray*}
\lim_{N \to \infty} \frac{1}{h} (LP_{\frac{\xi}{\sqrt{N}}, \alpha}[\widetilde{W}_2] - LP_{\frac{\xi}{\sqrt{N}}, \alpha}[W_2]) = 0
\end{eqnarray*}
where
$$LP_{\frac{\xi}{\sqrt{N}}, \alpha}[\widetilde{W}_2] = P_{\frac{\xi}{\sqrt{N}}}\left[ {V_{0}[\widetilde{W}_2]}^{-1/2}(\widetilde{W}_2 -E_0[\widetilde{W}_2]) > v_{1-\alpha}\right]$$
and $LP_{\frac{\xi}{\sqrt{N}}, \alpha}[W_2]$ is also the same.
\end{theorem}
\noindent{\it Proof.} See Appendices.

\section{Simulation study}
\label{sect:simul}

In this section, we compare the significance probabilities of $\widetilde{M}$ and $\widetilde{W}_2$ by simulation because the distributions of them depend on $F$. For 100,000 times random samples from the standard normal distribution, we estimate the significance probabilities in the tale area
$$\widetilde{\Omega}_{\alpha}=\left\{{\mbf x}\in{\mbf R}^n ~\left|~ \frac{\widetilde{m}({\mbf x})-E_1}{\sqrt{V_1}} > v_{1-\alpha},\right.\hspace{4mm}{\rm or}\hspace {5mm}\frac{\widetilde{w}_2({\mbf x})-E_2}{\sqrt{V_2}} > v_{1-\alpha}\right\}.$$
Similarly as Table 5, Table 6 shows the ratio of the samples that the significance probability of $\widetilde{W}_2$ is smaller than $\widetilde{M}$. Here we use the exponential type kernel as introduced in Remark 2.6 in $\widetilde{M}$, and the Epanechnikov kernel in $\widetilde{W}_2$. The both bandwidths are same, and $h=N^{-1/4}/ \log N$ as we will explain later. Comparing Table 5 and 6, we can see that the differences of the $p$-values of $\widetilde{M}$ and $\widetilde{W}_2$ is smaller than those of $M$ and $W_2$.

\begin{table}[h]
\caption{Relation of significance probabilities \label{tab:table5}}
\begin{center}
\begin{tabular}{c|cccccc}\hline
sample size ($m,n$)&(10,10)&(20,20)&(30,30)&(10,20)&(20,10)&($U_m^*,U_n^*$)\footnotemark[2]\\\hline
$z_{0.9}$ & 1.41 & 1.08 & 1.31 & 1.05 & 1.22 & 0.827\\
$z_{0.95}$ & 0.829 & 1.21 & 1.39 & 1.41 & 1.24 & 0.782\\
$z_{0.975}$ & 1.64 & 1.12 & 1.09 & 1.23 & 1.39 & 0.755\\
$z_{0.99}$ & 0.985 & 1.03 & 1.32 & 0.707 & 0.806 & 0.677\\\hline
\end{tabular}
\end{center}
\end{table}
\footnotetext[2]{$U_m^*$ and $U_n^*$ are random numbers from the discrete uniform distribution $U^*(5,40)$}

Next we will compare the powers of the smoothed sign, the smoothed Wilcoxon's and the two-sample $T$-test. The Pitman's efficiency of the $T$-test is given by

$$e_P[T_2] = \sqrt{\frac{\lambda(1-\lambda)}{V_0[X_1]}}.$$

In Table 7 and Table 8, we also use the exponential type kernel, the Epanechnikov kernel and $h=N^{-1/4}/ \log N$ which attains almost maximum order under the assumptions, to focus on the first orders of their local power. We simulate their power when $\theta= 0.1, ~0.5$ and the significance level $\alpha=0.01, ~0.05$, based on 100,000 repetitions.  In order to check the size condition, we also simulate the case $\theta=0$ in Table 9. When the underlying distribution $F$ is the heavy-tailed $T$-distributions, the simulation results show that the smoothed median test is superior than the other tests.  When the underlying distribution $F$ is the logistic, the simulation results show that the smoothed Wilcoxon's test is best. The student $t$-test is superior than others, when $F$ is normal. These simulation studies coincide with the Pitman's $A.R.E.$s, and justify the continuation of the ordinal tests. Also all of the empirical sizes of the proposed smoothed tests is close to those of the significance levels, and we can conclude they are asymptotically nonparametric.

\begin{remark}
One possible way of choosing their bandwidths is to reduce the errors of the approximations. It is especially important to `fix' the type I error in statistical testing, however there are no tradeoffs about the bandwidths in their normal approximations. Although the best $h$ is as small as possible from the asymptotic result, `smoothed' bootstrap method can give a numerical solution. Under the assumptions, there are no same values, and the appropriately smoothed resampling never returns `ties' similarly. If we choose the squared loss at the significance level $\alpha$, using $L$ repetition, the best bandwidth of $\widetilde{M}$ is given by 
$$\widetilde{h} = \min_{h>0} \left( \sharp\left\{ l : {V_{1}}^{-1/2}(\widetilde{M}_{\star}^{(l)} -E_1) > v_{1-\alpha} \right\} - L\alpha \right)^2$$
where $\widetilde{M}_{\star}^{(l)}$ is the $l$th test statistic given by $\{X_{\star,1}^{(l)}, \cdots, X_{\star,m}^{(l)}, X_{\star,m+1}^{(l)}, \cdots, X_{\star,m+n}^{(l)}\}$ and $X_{\star,1}^{(l)}, \cdots, X_{\star,m+n}^{(l)}$ is i.i.d. sample from smooth $\widetilde{F}$ estimated by the original $\{X_1, \cdots, X_m\}$. The minimizer is not decided uniquely because the function is discrete, so we choose the middle. If $L$ is enough large, the effect of choosing the middle is negligible.
\end{remark}

As we confirmed in Table 9, the normal approximation of $p$-values of $\widetilde{W}_2$ is good enough, but that of $\widetilde{M}$ is not. To see improvement of the size condition of $\widetilde{M}$, we do the numerical study, lastly. Using the standard normal density as the kernel function and cross-validated bandwidth, we obtain the resampled data, numerically optimal bandwidth and the results in Table 10 based on 100,000 repetitions. The bandwidth is calculated every time, using 1,000 resampled data sets. As we can see, the size conditions are improved and the numerical powers are also the same in these cases.

The kernel function does not affect much both theoretical and numerical result as same as the kernel density estimation under the assumption. The better choice of the bandwidths should be discussed more, but we postpone it for future work.

\newpage

\begin{table}[!h]
\caption{Power comparisons of $\widetilde{M}$, $\widetilde{W}_2$ and $T_2$ test \label{tab:table6}}
\begin{center}
\begin{tabular}{||c||c|c|c||c||c|c|c||}
   \hline
   $\alpha={0.01}$ & $m=30$ & $n$=30 &  &  &  &  &  \\
  \hline
  $\theta={0.1}$ & $\widetilde{M}$ & $\widetilde{W}_2$ & $T_2$ & $\theta={0.5}$ & $\widetilde{M}$ & $\widetilde{W}_2$ & $T_2$ \\
  \hline
  $N(0,1)$ & 0.01956 & 0.02564 & 0.02637 & $N(0,1)$ & 0.19579 & 0.31757 & 0.32987 \\
  Logis. & 0.01468 & 0.01794 & 0.01660 & Logis. & 0.07731 & 0.11056 & 0.10141  \\
  D.Exp. & 0.02230 & 0.02257 & 0.01892 & D.Exp. & 0.23392 & 0.23788 & 0.17014 \\
  \hline
\hline
   $\alpha=0.05$ & $m=30$ & $n$=30 &  &  &  &  &  \\
  \hline
  $\theta={0.1}$ & $\widetilde{M}$ & $\widetilde{W}_2$ & $T_2$ & $\theta={0.5}$ & $\widetilde{M}$ & $\widetilde{W}_2$ & $T_2$ \\
  \hline
  $N(0,1)$ & 0.08120 & 0.10458 & 0.10444 & $N(0,1)$ & 0.40059 & 0.59413 & 0.60390 \\
  Logis. & 0.06496 & 0.07868 & 0.07606 & Logis. & 0.20295 & 0.29871 & 0.28053 \\
  D.Exp. & 0.09006 & 0.09685 & 0.08602 & D.Exp. & 0.45217 & 0.49389 & 0.39608 \\
  \hline\hline    
   $\alpha={0.01}$ & $m=50$ & $n$=30 &  &  &  &  &  \\
  \hline
  $\theta={0.1}$ & $\widetilde{M}$ & $\widetilde{W}_2$ & $T_2$ & $\theta={0.5}$ & $\widetilde{M}$ & $\widetilde{W}_2$ & $T_2$ \\
  \hline
  $N(0,1)$ & 0.02294 & 0.02762 & 0.02926 & $N(0,1)$ & 0.24816 & 0.40303 & 0.41430 \\
  Logis. & 0.01703 & 0.01838 & 0.01820 & Logis. & 0.10060 & 0.13703 &  0.12568 \\
  D.Exp. & 0.02571 & 0.02495 & 0.02064 & D.Exp. & 0.29770 & 0.30226 & 0.21347 \\
  \hline
   $\alpha=0.05$ & $m=50$ & $n$=30 &  &  &  &  &  \\
  \hline
  $\theta={0.1}$ & $\widetilde{M}$ & $\widetilde{W}_2$ & $T_2$ & $\theta={0.5}$ & $\widetilde{M}$ & $\widetilde{W}_2$ & $T_2$ \\
  \hline
  $N(0,1)$ & 0.08854 & 0.11211 & 0.11319 & $N(0,1)$ & 0.48952 & 0.67834 & 0.68776 \\
  Logis. & 0.07280 & 0.08203 & 0.08026 & Logis. & 0.26507 & 0.34398 & 0.32542 \\
  D.Exp. & 0.09681 & 0.10313 & 0.09197 & D.Exp. & 0.54467 & 0.57204 & 0.45771 \\
    \hline\hline    
   $\alpha={0.01}$ & $m=30$ & $n$=50 &  &  &  &  &  \\
  \hline
  $\theta={0.1}$ & $\widetilde{M}$ & $\widetilde{W}_2$ & $T_2$ & $\theta={0.5}$ & $\widetilde{M}$ & $\widetilde{W}_2$ & $T_2$ \\
  \hline
  $N(0,1)$ & 0.02290 & 0.02793 & 0.02846 & $N(0,1)$ & 0.25171 & 0.40378 & 0.41573 \\
  Logis. & 0.01746 & 0.01856 & 0.01805 & Logis. & 0.09872 & 0.13701 & 0.12613 \\
  D.Exp. & 0.02692 & 0.02495 & 0.02094 & D.Exp. & 0.30090 & 0.30376 & 0.21466 \\
  \hline
   $\alpha=0.05$ & $m=30$ & $n$=50 &  &  &  &  &  \\
  \hline
  $\theta={0.1}$ & $\widetilde{M}$ & $\widetilde{W}_2$ & $T_2$ & $\theta={0.5}$ & $\widetilde{M}$ & $\widetilde{W}_2$ & $T_2$ \\
  \hline
  $N(0,1)$ & 0.09180 & 0.11229 & 0.11213 & $N(0,1)$ & 0.50837 & 0.67818 & 0.69049 \\
  Logis. & 0.07556 & 0.08381 & 0.07957 & Logis. & 0.27437 & 0.34509 & 0.32639 \\
  D.Exp. & 0.10089 & 0.10232 & 0.09172 & D.Exp. & 0.56812 & 0.57234 & 0.45807 \\
    \hline\hline    
   $\alpha={0.01}$ & $m=50$ & $n$=50 &  &  &  &  &  \\
  \hline
  $\theta={0.1}$ & $\widetilde{M}$ & $\widetilde{W}_2$ & $T_2$ & $\theta={0.5}$ & $\widetilde{M}$ & $\widetilde{W}_2$ & $T_2$ \\
  \hline
  $N(0,1)$ & 0.02659 & 0.03235 & 0.03267 & $N(0,1)$ & 0.35734 & 0.53628 & 0.55407 \\
  Logis. & 0.01974 & 0.02050 & 0.02001 & Logis. & 0.14500 & 0.18420 & 0.17026 \\
  D.Exp. & 0.03073 & 0.02898 & 0.02341 & D.Exp. & 0.43233 & 0.41255 & 0.28663 \\
  \hline
   $\alpha=0.05$ & $m=50$ & $n$=50 &  &  &  &  &  \\
  \hline
  $\theta={0.1}$ & $\widetilde{M}$ & $\widetilde{W}_2$ & $T_2$ & $\theta={0.5}$ & $\widetilde{M}$ & $\widetilde{W}_2$ & $T_2$ \\
  \hline
  $N(0,1)$ & 0.09766 & 0.12477 & 0.12720 & $N(0,1)$ & 0.60555 & 0.78490 & 0.79775 \\
  Logis. & 0.07692 & 0.08768 & 0.08545 & Logis. & 0.31511 & 0.41698 & 0.39451 \\
  D.Exp. & 0.11001 & 0.11382 & 0.09828 & D.Exp. & 0.67897 & 0.67842 & 0.54885 \\
\hline
\end{tabular}
\end{center}
\end{table}

\begin{table}[!h]
\caption{Power comparisons of $\widetilde{M}$, $\widetilde{W}_2$ and $T_2$ test \label{tab:table6}}
\begin{center}
\begin{tabular}{||c||c|c|c||c||c|c|c||}
   \hline
   $\alpha=0.01$ & $m=30$ & $n$=30 &  &  &  &  &  \\
  \hline
  $\theta={0.1}$ & $\widetilde{M}$ & $\widetilde{W}_2$ & $T_2$ & $\theta={0.5}$ & $\widetilde{M}$ & $\widetilde{W}_2$ & $T_2$ \\
  \hline
  $T(2)$ & 0.01596 & 0.01953 & 0.01177 & $T(2)$ & 0.14133 & 0.16649 & 0.07647 \\
  $T(1)$ & 0.01533 & 0.01835 & 0.00314 & $T(1)$ & 0.11032 & 0.10041 & 0.01019 \\
  $T(1/2)$ & 0.01365 & 0.01466 & 0.00025 & $T(1/2)$ & 0.07537 & 0.04999 & 0.00037 \\
\hline
   $\alpha=0.05$ & $m=30$ & $n$=30 &  &  &  &  &  \\
  \hline
  $\theta={0.1}$ & $\widetilde{M}$ & $\widetilde{W}_2$ & $T_2$ & $\theta={0.5}$ & $\widetilde{M}$ & $\widetilde{W}_2$ & $T_2$ \\
  \hline
  $T(2)$ & 0.06720 & 0.08651 & 0.06764 & $T(2)$ & 0.32518 & 0.39554 & 0.23547 \\
  $T(1)$ & 0.06375 & 0.07879 & 0.03502 & $T(1)$ & 0.27208 & 0.27884 & 0.06869 \\
  $T(1/2)$ & 0.05613 & 0.06889 & 0.01277 & $T(1/2)$ & 0.19935 & 0.16884 & 0.01445 \\
  \hline\hline    
   $\alpha={0.01}$ & $m=50$ & $n$=30 &  &  &  &  &  \\
  \hline
  $\theta={0.1}$ & $\widetilde{M}$ & $\widetilde{W}_2$ & $T_2$ & $\theta={0.5}$ & $\widetilde{M}$ & $\widetilde{W}_2$ & $T_2$ \\
  \hline
  $T(2)$ & 0.02088 & 0.02185 & 0.01263 & $T(2)$ & 0.19386 & 0.21044 & 0.09175 \\
  $T(1)$ & 0.02009 & 0.01878 & 0.00312 & $T(1)$ & 0.15396 & 0.12362 & 0.01078 \\
  $T(1/2)$ & 0.01828 & 0.01498 & 0.00020 & $T(1/2)$ & 0.10638 & 0.05750 & 0.00028 \\
  \hline
   $\alpha=0.05$ & $m=50$ & $n$=30 &  &  &  &  &  \\
  \hline
  $\theta={0.1}$ & $\widetilde{M}$ & $\widetilde{W}_2$ & $T_2$ & $\theta={0.5}$ & $\widetilde{M}$ & $\widetilde{W}_2$ & $T_2$ \\
  \hline
  $T(2)$ & 0.07719 & 0.09306 & 0.07151 & $T(2)$ & 0.40557 & 0.45878 & 0.26322 \\
  $T(1)$ & 0.07433 & 0.08112 & 0.03657 & $T(1)$ & 0.34600 & 0.32085 &  0.06955 \\
  $T(1/2)$ & 0.07085 & 0.07098 & 0.01368 & $T(1/2)$ & 0.26580 & 0.18676 & 0.01427 \\
    \hline\hline    
   $\alpha={0.01}$ & $m=30$ & $n$=50 &  &  &  &  &  \\
  \hline
  $\theta={0.1}$ & $\widetilde{M}$ & $\widetilde{W}_2$ & $T_2$ & $\theta={0.5}$ & $\widetilde{M}$ & $\widetilde{W}_2$ & $T_2$ \\
  \hline
  $T(2)$ & 0.02242 & 0.02156 & 0.01249 & $T(2)$ & 0.20139 & 0.21035 & 0.09152 \\
  $T(1)$ & 0.02166 & 0.01941 & 0.00318 & $T(1)$ & 0.15691 & 0.12394 & 0.01026 \\
  $T(1/2)$ & 0.01924 & 0.01513 & 0.00031 & $T(1/2)$ & 0.10540 & 0.05795 & 0.00031 \\
  \hline
   $\alpha=0.05$ & $m=30$ & $n$=50 &  &  &  &  &  \\
  \hline
  $\theta={0.1}$ & $\widetilde{M}$ & $\widetilde{W}_2$ & $T_2$ & $\theta={0.5}$ & $\widetilde{M}$ & $\widetilde{W}_2$ & $T_2$ \\
  \hline
  $T(2)$ & 0.08220 & 0.09153 & 0.07086 & $T(2)$ & 0.42438 & 0.45922 & 0.26436 \\
  $T(1)$ & 0.07868 & 0.08164 & 0.03563 & $T(1)$ & 0.36214 & 0.32053 & 0.06922 \\
  $T(1/2)$ & 0.07288 & 0.07051 & 0.01360 & $T(1/2)$ & 0.27672 & 0.18796 & 0.01486 \\
    \hline\hline    
   $\alpha={0.01}$ & $m=50$ & $n$=50 &  &  &  &  &  \\
  \hline
  $\theta={0.1}$ & $\widetilde{M}$ & $\widetilde{W}_2$ & $T_2$ & $\theta={0.5}$ & $\widetilde{M}$ & $\widetilde{W}_2$ & $T_2$ \\
  \hline
  $T(2)$ & 0.02325 & 0.02491 & 0.01386 & $T(2)$ & 0.26932 & 0.28553 & 0.11333 \\
  $T(1)$ & 0.02127 & 0.02060 & 0.00313 & $T(1)$ & 0.21263 & 0.16481 &  0.01071 \\
  $T(1/2)$ & 0.02020 & 0.01649 & 0.00033 & $T(1/2)$ & 0.14431 & 0.07195 & 0.00024 \\
  \hline
   $\alpha=0.05$ & $m=50$ & $n$=50 &  &  &  &  &  \\
  \hline
  $\theta={0.1}$ & $\widetilde{M}$ & $\widetilde{W}_2$ & $T_2$ & $\theta={0.5}$ & $\widetilde{M}$ & $\widetilde{W}_2$ & $T_2$ \\
  \hline
  $T(2)$ & 0.08884 & 0.10073 & 0.07631 & $T(2)$ & 0.50420 & 0.55280 & 0.30033 \\
  $T(1)$ & 0.08219 & 0.08522 & 0.03600 & $T(1)$ & 0.42768 & 0.38690 &  0.07028 \\
  $T(1/2)$ & 0.07409 & 0.07254 & 0.01420 & $T(1/2)$ & 0.32194 & 0.22081 & 0.01461 \\
\hline
\end{tabular}
\end{center}
\end{table}

\begin{table}[!h]
\caption{Check on the size condition of $\widetilde{M}$, $\widetilde{W}_2$ and $T_2$ test \label{tab:table6}}
\begin{center}
\begin{tabular}{||c||c|c|c||c||c|c|c||}
   \hline
   $\alpha={0.01}$ & $m=30$ & $n$=30 &  &  &  &  &  \\
  \hline
  $\theta={0}$ & $\widetilde{M}$ & $\widetilde{W}_2$ & $T_2$ & $\theta={0}$ & $\widetilde{M}$ & $\widetilde{W}_2$ & $T_2$ \\
  \hline
  $N(0,1)$ & 0.00955 & 0.01006 & 0.01010 & $T(2)$ & 0.00776 & 0.01028 & 0.00651 \\
  Logis. & 0.00873 & 0.01015 & 0.00946 & $T(1)$ & 0.00765 & 0.00999 &  0.00231 \\
  D.Exp. & 0.00921 & 0.01002 & 0.00859 & $T(1/2)$ & 0.00780 & 0.00930 & 0.00014 \\
  \hline
\hline
   $\alpha=0.05$ & $m=30$ & $n$=30 &  &  &  &  &  \\
  \hline
  $\theta={0}$ & $\widetilde{M}$ & $\widetilde{W}_2$ & $T_2$ & $\theta={0}$ & $\widetilde{M}$ & $\widetilde{W}_2$ & $T_2$ \\
  \hline
  $N(0,1)$ & 0.04647 & 0.05135 & 0.04854 & $T(2)$ & 0.03893 & 0.05124 & 0.04530 \\
  Logis. & 0.04470 & 0.05115 & 0.04922 & $T(1)$ & 0.03791 & 0.05087 & 0.02955 \\
  D.Exp. & 0.04609 & 0.05153 & 0.04928 & $T(1/2)$ & 0.03765 & 0.05086 & 0.01282 \\
  \hline\hline    
   $\alpha={0.01}$ & $m=50$ & $n$=30 &  &  &  &  &  \\
  \hline
  $\theta={0}$ & $\widetilde{M}$ & $\widetilde{W}_2$ & $T_2$ & $\theta={0}$ & $\widetilde{M}$ & $\widetilde{W}_2$ & $T_2$ \\
  \hline
  $N(0,1)$ & 0.01003 & 0.01029 & 0.01000 & $T(2)$ & 0.01034 & 0.01061 & 0.00660 \\
  Logis. & 0.01023 & 0.01010 & 0.00960 & $T(1)$ & 0.01021 & 0.00981 & 0.00219 \\
  D.Exp. & 0.00960 & 0.01008 & 0.00846 & $T(1/2)$ & 0.01061 & 0.00942 & 0.00020 \\
  \hline
   $\alpha=0.05$ & $m=50$ & $n$=30 &  &  &  &  &  \\
  \hline
  $\theta={0}$ & $\widetilde{M}$ & $\widetilde{W}_2$ & $T_2$ & $\theta={0}$ & $\widetilde{M}$ & $\widetilde{W}_2$ & $T_2$ \\
  \hline
  $N(0,1)$ & 0.04567 & 0.05033 & 0.04961 & $T(2)$ & 0.04248 & 0.05108 & 0.04512 \\
  Logis. & 0.04742 & 0.05056 & 0.04980 & $T(1)$ & 0.04248 & 0.05032 & 0.02934 \\
  D.Exp. & 0.04496 & 0.05169 & 0.04836 & $T(1/2)$ & 0.04452 & 0.05000 & 0.01291 \\
    \hline\hline
   $\alpha={0.01}$ & $m=30$ & $n$=50 &  &  &  &  &  \\
  \hline
  $\theta={0}$ & $\widetilde{M}$ & $\widetilde{W}_2$ & $T_2$ & $\theta={0}$ & $\widetilde{M}$ & $\widetilde{W}_2$ & $T_2$ \\
  \hline
  $N(0,1)$ & 0.01031 & 0.00985 & 0.01002 & $T(2)$ & 0.01068 & 0.01039 & 0.00668 \\
  Logis. & 0.00985 & 0.01010 & 0.00949 & $T(1)$ & 0.01028 & 0.00925 & 0.000246 \\
  D.Exp. & 0.01012 & 0.00983 & 0.00905 & $T(1/2)$ & 0.01052 & 0.00990 & 0.00030 \\
  \hline
   $\alpha=0.05$ & $m=30$ & $n$=50 &  &  &  &  &  \\
  \hline
  $\theta={0}$ & $\widetilde{M}$ & $\widetilde{W}_2$ & $T_2$ & $\theta={0}$ & $\widetilde{M}$ & $\widetilde{W}_2$ & $T_2$ \\
  \hline
  $N(0,1)$ & 0.04760 & 0.05032 & 0.05034 & $T(2)$ & 0.04514 & 0.05126 & 0.04625 \\
  Logis. & 0.04943 & 0.05061 & 0.05003 & $T(1)$ & 0.04553 & 0.05031 & 0.03082 \\
  D.Exp. & 0.04727 & 0.05196 & 0.04883 & $T(1/2)$ & 0.04664 & 0.05043 & 0.01316 \\
    \hline\hline    
   $\alpha={0.01}$ & $m=50$ & $n$=50 &  &  &  &  &  \\
  \hline
  $\theta={0}$ & $\widetilde{M}$ & $\widetilde{W}_2$ & $T_2$ & $\theta={0}$ & $\widetilde{M}$ & $\widetilde{W}_2$ & $T_2$ \\
  \hline
  $N(0,1)$ & 0.01061 & 0.00995 & 0.01035 & $T(2)$ & 0.00962 & 0.01040 & 0.00698 \\
  Logis. & 0.01117 & 0.01006 & 0.00990 & $T(1)$ & 0.00947 & 0.00959 & 0.00237 \\
  D.Exp. & 0.00988 & 0.01029 & 0.00919 & $T(1/2)$ & 0.01020 & 0.00984 & 0.00032 \\
  \hline
   $\alpha=0.05$ & $m=50$ & $n$=50 &  &  &  &  &  \\
  \hline
  $\theta={0}$ & $\widetilde{M}$ & $\widetilde{W}_2$ & $T_2$ & $\theta={0}$ & $\widetilde{M}$ & $\widetilde{W}_2$ & $T_2$ \\
  \hline
  $N(0,1)$ & 0.04696 & 0.05005 & 0.05066 & $T(2)$ & 0.04661 & 0.05080 & 0.04639 \\
  Logis. & 0.04707 & 0.05061 & 0.05010 & $T(1)$ & 0.04506 & 0.05084 & 0.03048 \\
  D.Exp. & 0.04650 & 0.05212 & 0.04815 & $T(1/2)$ & 0.04538 & 0.05061 & 0.01344 \\\hline
\end{tabular}
\end{center}
\end{table}

\clearpage

\begin{table}[!h]
\caption{Numerical power of $\widetilde{M}$ using asymptotically optimal  bandwidth \label{tab:table6}}
\begin{center}
\begin{tabular}{||c||c|c|c||c|c|c||}
  \hline
  $\alpha=0.05$ & $\theta={0}$ & $\theta={0.1}$ & $\theta={0.5}$ & $\theta={0}$ & $\theta={0.1}$ & $\theta={0.5}$ \\
   \hline
  $\widetilde{M}$ & $m=30$ & $n$=30 &  & $m=50$ & $n=30$ & \\
  \hline
  $T(2)$ & 0.05292 & 0.08383 & 0.35837 & 0.05004 & 0.08798 & 0.41117 \\
  $T(1)$ & 0.05325 & 0.08237 & 0.31159 & 0.05160 & 0.08571 & 0.35859 \\
  $T(1/2)$ & 0.04695 & 0.06989 & 0.24244 & 0.05147 & 0.07699 & 0.28194 \\
\hline
  $\widetilde{M}$ & $m=30$ & $n$=50 &  & $m=50$ & $n=50$ & \\
  \hline
  $T(2)$ & 0.04871 & 0.08919 & 0.42320 & 0.04901 & 0.09598 & 0.50938 \\
  $T(1)$ & 0.05220 & 0.08731 & 0.36653 & 0.05250 & 0.09441 & 0.44241 \\
  $T(1/2)$ & 0.04824 & 0.07441 & 0.28882 & 0.04932 & 0.08373 & 0.35290 \\
\hline
\end{tabular}
\end{center}
\end{table}

\bibliographystyle{gSTA}
\bibliography{ref}

\section{Appendices}

{\bf Proof of Lemma 2.2}

By the direct computation, we can find
\begin{eqnarray*}
&&F_{Z}^{O}(z) -F_{Z}^{E}(z) \\
&=& \frac{1}{\beta(N/2,N/2)} \int_{-\infty}^z [F(x)]^{(N/2)-1} f(x) \biggm\{ \frac{\beta(N/2,N/2)}{\beta((N+1)/2,(N+1)/2)} \\
&& ~~~~~~ \Bigm(\sqrt{F(x)} [1-F(x)]^{(N-1)/2} \Bigm) -2\Bigm( [1-F(x)]^{(N/2)}  [1-F(2z -x)]^{(N/2)} \Bigm) \biggm \} dx.
\end{eqnarray*}
Since
$$\lim_{N \rightarrow \infty} \frac{\Gamma(N + \eta)}{\Gamma(N) N^{\eta}} =1 ~~~ (\eta \in {\mathbf R}),$$
we can obtain
$$\lim_{N \rightarrow \infty} \frac{\beta(N/2,N/2)}{\beta((N+1)/2,(N+1)/2)} =1.$$
When $[1-F(x)]^{N} \not\to 0$, it is easily to see $F^N(x) \to 0$ and vice versa. If both fails, both $[F(x)]^{N}$ and $[1-F(x)]^{N}$ shrinks to $0$. From the above, we can see the result.

{\bf Proof of Theorem 2.3} 

Hereafter, we assume that $h = o(N^{-1/2})$ or that $A_{1,1}^*=A_{1,2}^*=0$ and $h = o(N^{-1/3})$, to ignore the residual term.

We utilize the following Bahadur representation (Bahadur \cite{bahadur1966note}) of the combined median $Z$ under $H_0$
$$Z = z_0 + \frac{1}{f(z_0)} \left[\frac{1}{2} -F_{(N)}(z_0) \right] + R_N$$
where $z_0$ is the median of $F$, $F_{(N)}$ stands for the empirical distribution function of $\{Z_1, \cdots ,Z_N\}$ and $R_N$ is the residual which satisfies $R_N=o_P(N^{-3/4} \log N)$.
By the following moment evaluation of the residual term of $Z$ (Reiss \cite{reiss1976asymptotic})
$$N^{l/2}E[(Z-z_0)^l] = (p(1-p))^{l/2} \frac{\kappa_{l}}{f^{l}(z_0)} + O(N^{-1/2})$$
and
$$E[Z-z_0] = -\frac{f^{(1)}(z_0)}{f^3(z_0)} E[(B-z_0)^2] + o(N^{-1})$$
where $\kappa_{l}$ is the $l$-$th$ moment of the standard normal distribution and $B$ stands for the beta distribution $Beta(N/2,N/2)$, we find
$$E[R_N] = -\frac{f^{(1)}(z_0)}{4N f^3(z_0)} + O(N^{-2})$$
and
$$E[R_N^2]  = o(N^{-1}).$$

Hereafter we put $N$ is odd. One of $\{X_i\}_{i=1, \cdots m}$ may be $Z$ but we don't need the differentiability of $K^*$. This is because the probability of $Z-X_i$ takes any finite points ($\neq 0$) is 0, and to smooth $K^*$ near $0$ and $1$ appropriately is possible. Using the representation, we can see the following asymptotic expansion under $H_0$
\begin{eqnarray*}
\widetilde{M} &=&\sum_{i=1}^m \left[ K^*\left(\frac{z_0 -X_i}{h}\right) +\frac{1}{h}k^*\left(\frac{z_0 -X_i}{h}\right) (Z -z_0) + \frac{1}{h^2} {k^*}^{(1)} \left(\frac{z_0 -X_i}{h}\right) (Z -z_0)^2 \right] \\
&& ~~~~~~ + o_P(1) \\
&=&\sum_{i=1}^m \biggm[ K^*\left(\frac{z_0 -X_i}{h}\right) +\frac{1}{h}k^*\left(\frac{z_0 -X_i}{h}\right) \left[ \frac{1}{f(z_0)} \left( \frac{1}{2} -F_{(N)}(z_0) \right) + R_N \right] \\
&& ~~~~~~ + \frac{1}{f^2(z_0) h^2} {k^*}^{(1)} \left(\frac{z_0 -X_i}{h}\right) \left[\frac{1}{2} -F_{(N)}(z_0) \right]^2 \biggm] + o_P(1)
\end{eqnarray*}
where ${k^*}^{(1)}$ is the derivative of ${k^*}$, $F_{(N)}$ is the empirical distribution function of $\{Z_1, \cdots, Z_N\}$ and $R_N$ is a residual which satisfies $R_N = o_P(N^{-3/4} \log N)$.

Using conditional expectation given $X_j$ and the Taylor expansion, we obtain
\begin{eqnarray*}
E_0[\widetilde{M}] &=& m \int_{-\infty}^{z_0} \biggm[K^*\left(\frac{z_0 -x}{h}\right) +\frac{1}{h}k^*\left(\frac{z_0 -x}{h}\right) \biggm\{ \frac{1}{f(z_0) N}\left(\frac{1}{2} -\psi(z_0-x) \right) \\
&& ~~~~~~~~~~~~ + E[R_N]\biggm \} \biggm]f(x)dx \\
&& ~~~~~~ + \frac{m}{f^2(z_0) N h^2} E\left[ {k^*}^{(1)} \left(\frac{z_0 -X_i}{h}\right) \left(\frac{1}{2} -\psi(X_j-z_0) \right)^2 \right]_{i \neq j} \\
&& ~~~~~~ + o(1 +\log N (N^{1/4}h^2 +N^{-1/4} h^{-1/2}))\\
&=& m \int_{0}^{\infty} k^*\left(v\right) \biggm[F(z_0 -hv)+ f(z_0 -hv) \left\{-\frac{1}{2f(z_0) N} -\frac{f^{(1)}(z_0)}{4N f^3(z_0)} \right\} \\
&& ~~~~~~ + \frac{m}{4N f^2(z_0)} f^{(1)}(z_0 -hv)\biggm] dv  + o(1 +(\log N) (N^{1/4}h^2 +N^{-1/4} h^{-1/2}))\\
&=& \frac{m}{2} \left[1 + h^2 A_{2,1}^* f^{(1)}(z_0) - \frac{1}{m+n} \right]  + o(1 + (\log N) N^{-1/4} h^{-1/2}).
\end{eqnarray*}

Using the stochastic expansion and the above results, we have
\begin{eqnarray*}
\widetilde{M}^2 &=& \Biggm[ \sum_{i=1}^m \biggm\{ K^*\left(\frac{z_0 -X_i}{h}\right) +\frac{1}{h}k^*\left(\frac{z_0 -X_i}{h}\right) \left[ \frac{1}{f(z_0)} \left( \frac{1}{2} -F_{(N)}(z_0) \right) + R_N \right] \\
&& + \frac{1}{f^2(z_0) h^2} {k^*}^{(1)} \left(\frac{z_0 -X_i}{h}\right)  \left[\frac{1}{2} -F_{(N)}(z_0) \right]^2 \biggm\} + O_P(\sqrt{N} R_N + N R_N^2 +N^{-1/2}) \biggm]^2.
\end{eqnarray*}

From the Bahadur representation, we can obtain
\begin{eqnarray*}
&&E \left[ \sum_{i=1}^m K^*\left(\frac{z_0 -X_i}{h}\right) O_P(\sqrt{N} R_N + N R_N^2 +N^{-1/2}) \right] \\
&=& O(N^{1/2} \log N),
\end{eqnarray*}
so the expectation of the squared value is given by
\begin{eqnarray*}
&& E_0[\widetilde{M}^2] \\
&=& \sum_{i=1}^m \sum_{j=1}^m \Biggm( E\left[ K^*\left(\frac{z_0 -X_i}{h}\right) K^*\left(\frac{z_0 -X_j}{h}\right) \right] \\
&&+\frac{2}{h} E\left[K^*\left(\frac{z_0 -X_i}{h}\right) k^*\left(\frac{z_0 -X_j}{h}\right) \left( \frac{1}{f(z_0) N} (1 -\psi(z_0-X_i) -\psi(z_0-X_j)) +R_N \right) \right] \\
&&+ \frac{1}{f^2(z_0) (Nh)^2} \sum_{k=1}^N \sum_{l=1}^N \\
&& ~~~ \Biggm\{ E\left[K^*\left(\frac{z_0 -X_i}{h}\right) {k^*}^{(1)}\left(\frac{z_0 -X_j}{h}\right) \left\{\frac{1}{2} -\psi(z_0-Z_k)\right\} \left\{\frac{1}{2} -\psi(z_0-Z_l)\right\} \right] \\
&& ~~~~~ + 2 E\left[k^*\left(\frac{z_0 -X_i}{h}\right) k^*\left(\frac{z_0 -X_j}{h}\right) \left\{\frac{1}{2} -\psi(z_0-Z_k)\right\} \left\{\frac{1}{2} -\psi(z_0-Z_l)\right\} \right] \Biggm\} \Biggm)\\
&& \hspace{-2ex}+ O(N^{1/2} \log N).
\end{eqnarray*}

Calculating the expectations, we have the variance as follows
\begin{eqnarray*}
V_0[\widetilde{M}] &=& E_0[\widetilde{M}^2] - E_0[\widetilde{M}]^2 \\
&=& (m^2 -m) E\left[ {K^*}\left(\frac{z_0 -X_i}{h}\right) \right]^2 +m E\left[ {K^*}^2\left(\frac{z_0 -X_i}{h}\right) \right] \\
&& +\frac{2 m^2}{f(z_0) Nh} E\left[K^*\left(\frac{z_0 -X_i}{h}\right) k^*\left(\frac{z_0 -X_j}{h}\right) \left\{ 1 -\psi(z_0-X_i) -\psi(z_0-X_j) \right\} \right]_{i \neq j} \\
&&+\frac{2 m^2}{h} E\left[K^*\left(\frac{z_0 -X_i}{h}\right) k^*\left(\frac{z_0 -X_j}{h}\right) R_N \right]_{i \neq j} \\
&&+ \frac{2 m^2 N}{f^2(z_0) (Nh)^2} E\left[K^*\left(\frac{z_0 -X_i}{h}\right) {k^*}^{(1)} \left(\frac{z_0 -X_j}{h}\right) \left\{\frac{1}{2} -\psi(z_0-Z_k)\right\}^2 \right]_{i \neq j}\\
&&+ \frac{m^2 N}{f^2(z_0) (Nh)^2} E\left[k^*\left(\frac{z_0 -X_i}{h}\right) k^*\left(\frac{z_0 -X_j}{h}\right) \left\{\frac{1}{2} -\psi(z_0-Z_k)\right\}^2 \right]_{i \neq j} \\
&&-\frac{m^2}{4} \left[1 +h^2 A_{2,1}^* f^{(1)}(z_0) - \frac{1}{m+n} \right]^2 \\
&&+ O(Nh^2 +N^{1/2} \log N) \\
&=& \frac{m}{4} +\frac{m^2}{2N} -\frac{m^2}{N} -\frac{m^2 f^{(1)}(z_0)}{4 Nf^2(z_0)} +\frac{m^2 f^{(1)}(z_0)}{4 Nf^2(z_0)} +\frac{m^2}{4N} + O(Nh +N^{1/2} \log N) \\
&=& \frac{m}{4}- \frac{m^2}{4N} + O(Nh +N^{1/2} \log N) = \frac{mn}{4(m+n)} + O(Nh +N^{1/2} \log N).
\end{eqnarray*}

Using the Bahadur representation for quantiles of two samples  (Liu \& Yin \cite{liu1994asymptotic}) under $H_1$, and both $f(z_{\theta,N})$ and $f(z_{\theta,N} -\theta)$ are positive, we have
\begin{eqnarray*}
\widetilde{M} &=&\sum_{i=1}^m \left[ K^*\left(\frac{z_{\theta,N} -X_i}{h}\right) +\frac{1}{g_{\theta,N}(z_{\theta,N}) h}k^*\left(\frac{z_{\theta,N} -X_i}{h}\right) \left\{\frac{1}{2} -G_{\theta,(N)}(z_{\theta,N})\right\} \right] \\
&&+ o_P(N^{1/2})
\end{eqnarray*}
where
$$G_{\theta,(N)}(x) = \lambda_N F_{X,(m)}(x) +(1-\lambda_N) F_{Y,(n)}(x)$$
and $F_{X,(m)}$, $F_{Y,(n)}$ are the empirical distribution function of $\{X_1, \cdots ,X_m\}$ and $\{Y_1, \cdots ,Y_n\}$ respectively. Therefore, we see the following asymptotic expectation under $H_1$
\begin{eqnarray*}
E_{\theta}[\widetilde{M}] &=& m \int_{-\infty}^{\infty} k^*\left(v\right) F(z_{\theta,N} -hv)dy  + o(N^{1/2})\\
&=& m F(z_{\theta,N}) +o(N^{1/2}).
\end{eqnarray*}
Combining the above results, we can prove that the Pitman efficiency is same as the discrete one.

{\bf Proof of Theorem 2.4}

Using the result of Lahiri \& Chatterjee \cite{lahiri2007berry}, we have the following Berry-Esseen bound
$$\sup_{-\infty < x < \infty} \left|P_0 \left[ V_{3}^{-1/2} (M^{\dagger} -E_3) < x\right] -\Phi(x) \right| = O \left(\frac{1}{\sqrt{N}} \right) $$
where $E_3$ and $V_3$ are the expectation and the variance of $M^{\dagger}$. From the calculation of the proof of Theorem 2.3, we have
$$E_0[ D^2 ] = O\left( \sum_{i,j=1} (A_{i,1}^* h^i + A_{j,1,1}^* h^j) \right)$$
where $D=V_{1}^{-1/2} (\widetilde{M} -E_0[\widetilde{M}]) -V_{3}^{-1/2} (M^{\dagger} -E_3)$. Therefore, we can obtain $D = O_P(\sqrt{h})$ and $P[|D|>N^{-\epsilon} h^{1/4}] = N^{\epsilon} h^{1/4}$. Combining
$$\sup_{|t|<N^{-\epsilon} h^{1/2}, -\infty < x < \infty} \left| \Phi(x+t)-\Phi(x) \right| = O (N^{-\epsilon} h^{1/4}) $$
and the above result, we get 
\begin{eqnarray*}
&&\sup_{-\infty < x < \infty} \left|P_0 \left[ V_{1}^{-1/2} (\widetilde{M} -E_1) < x\right] -\Phi(x) \right| \\
&=&\sup_{-\infty < x < \infty} \left|P_0 \left[ \{V_{3}^{-1/2} (M^{\dagger} -E_3) + D + O(\sqrt{N} \sum_{i=1} A_{i,1}^* h^i)\} < x\right] -\Phi(x) \right| \\
&=& O \left(\sum_{i=1}(\sqrt{N} A_{i,1}^* h^i) + N^{\epsilon} h^{1/4} +\frac{1}{\sqrt{N}} \right).
\end{eqnarray*}
This completes the proof.

{\bf Proof of Theorem 2.5}

Put $\theta = N^{-1/2}\xi$ and $z_{\theta,N} = z_0 + N^{-1/2} \eta + o(N^{-1/2} h)$ where $\eta$ is a constant number and $z_0$ is defined before. Since the following expansion holds
\begin{eqnarray*}
\frac{N}{2} &=& m F(z_{\theta,N}) + n F(z_{\theta,N} -\theta)\\
&=& N F(z_{\theta,N}) - \frac{n}{\sqrt{N}} f(z_{\theta,N}) \xi + O(N^{-1}),
\end{eqnarray*}
we have
$$\sqrt{N} f(z_0) \eta = \frac{n}{\sqrt{N}} f(z_0) \xi + o(N^{1/2} h).$$
Then, we find $z_{\theta,N} = z_0 + N^{-1/2} (1- \lambda) \xi + o(N^{-1/2} h).$

Combining
$$V_0[\widetilde{M}] = \frac{mn}{4(m+n)} - 2(mh) A_{1,1,1}^* f(z_0) + O\left(Nh^2 +\frac{\log N}{\sqrt{N}} \right)$$
where $A_{i,j,l}^* = \int t^i {k^{*}}^j(t) {K^{*}}^l(t) dt$, and
\begin{eqnarray*}
E_{\theta}[\widetilde{M}] -E_{0}[\widetilde{M}] &=& m F(z_{\theta,N}) - \frac{m}{2} + O(N h^2)\\
&=& \lambda (1 -\lambda) \sqrt{N} f(z_0) \xi + o(N^{1/2} h),\\
\end{eqnarray*}
we can find that the local power of $\widetilde{M}$ is given by
$$LP_{\frac{\xi}{\sqrt{N}}, \alpha}[\widetilde{M}] = P_{\frac{\xi}{\sqrt{N}}}\left[ {V_{1}^{\star}}^{-1/2}(\widetilde{M} -E_1) > v_{1-\alpha}\right] = 1- \Phi( v_{1-\alpha} -c) + O(N^{-1/2} + h^2)$$
where
$$V_{1}^{\star} = V_{1} - 2(mh) A_{1,1,1}^* f(z_0),$$
$$c = \xi \left[ 2\sqrt{\lambda (1- \lambda)} f(z_0) + 8h \sqrt{\frac{\lambda}{1 -\lambda}} A_{1,1,1}^* f^2(z_0) \right]$$
and $v_{1-\alpha}$ is the $(1-\alpha)th$ quantile of the standard normal distribution. $M^{\dagger}$ doesn't have the second term of $c$ of order $h$, so we can find the result.

{\bf Proof of Theorem 3.1, 3.2 and 3.3}

We can easily see that the main term of the variance vanishes, and
\begin{eqnarray*}
V_{\theta}[\widetilde{W}_2] &=& \sum_{i=1}^m \sum_{k=1}^m \sum_{j=1}^n \sum_{l=1}^n E\left[ K\left(\frac{X_i-Y_j}{h}\right) K\left(\frac{X_k-Y_l}{h}\right) \right] - \left[mn \int_{-\infty}^{\infty} f(y) F(y-\theta) dy \right]^2 \\
&& ~~~ +o(N^3)\\
&=& m^2 n Cov_{\theta}\left( K\left(\frac{X_i-Y_j}{h}\right), K\left(\frac{X_k-Y_j}{h}\right) \right) \\
&& ~~~ + m n^2 Cov_{\theta}\left( K\left(\frac{X_i-Y_j}{h}\right), K\left(\frac{X_i-Y_l}{h}\right) \right) +o(N^3).
\end{eqnarray*}

By the direct calculation, we have
\begin{eqnarray*}
&&E_{\theta}\left[ K\left(\frac{X_i-Y_j}{h}\right) K\left(\frac{X_i-Y_l}{h}\right) \right] \\
&=& \int \cdots \int K\left(\frac{x -y}{h}\right) K\left(\frac{x -w}{h}\right) f(x) f(y-\theta) f(w-\theta)dx dy dw \\
&=& \int_{-\infty}^{\infty}  f(x) F(x-\theta) F(x-\theta)dx + O(h^2)
\end{eqnarray*}
and
\begin{eqnarray*}
&&E_{\theta}\left[ K\left(\frac{X_i-Y_j}{h}\right) K\left(\frac{X_k-Y_j}{h}\right) \right] \\
&=& \int \cdots \int K\left(\frac{x -y}{h}\right) K\left(\frac{u -y}{h}\right) f(x) f(u) f(y-\theta) dx du dy \\
&=& \int_{-\infty}^{\infty} \{1-F(y)\}^2 f(y-\theta) dy + O(h^2).
\end{eqnarray*}

Now, we can easily obtain the following equation
\begin{eqnarray*}
V_{\theta}[\widetilde{W}_2] &=& m^2 n \left( \int_{-\infty}^{\infty} f(x) F(x-\theta) F(x-\theta)dx - \left[ \int_{-\infty}^{\infty} f(y) F(y-\theta) dy \right]^2 \right) \\
&&+ m n^2 \left( \int_{-\infty}^{\infty} \{1-F(y)\}^2 f(y-\theta) dy- \left[\int_{-\infty}^{\infty} f(y) F(y-\theta) dy \right]^2 \right) +o(N^3).
\end{eqnarray*}

From the above, Theorem 3.1 is proved. Next, using the result of Maesono \cite{maesono1985edgeworth} and Garc\'{i}a-Soid\'{a}n et al \cite{garcia1997edgeworth}, we can obtain
$$\sup_{-\infty < x < \infty} \left|P_0 \left[ V_2^{-1/2} (\widetilde{W}_2 -E_2) < x\right] -Q_{m,n}(x) \right| = o(N^{-1/2})$$
where
\begin{eqnarray*}
Q_{m,n}(x) &=& \Phi(x) -\phi(x) \frac{x^2 -1}{6 \tau_{m,n}^3} \biggm[\frac{1}{m^2}E_0[g_{1,0}^3 (X)] + \frac{1}{n^2} E_0[g_{0,1}^3 (Y)] \\
&& ~~~~~~ + \frac{6}{mn} E_0[g_{1,0} (X) g_{0,1} (Y) g_{1,1} (X,Y)] \biggm],
\end{eqnarray*}
$$g_{1,0}(X) = \int k(v) F(X +hv) dv -\frac{1}{2} + O(C(h)),$$
$$g_{0,1}(Y) = \int k(v) F(Y -hv) dv -\frac{1}{2} + O(C(h)),$$
$$g_{1,1}(X,Y) = W\left(\frac{Y-X}{h}\right) -g_{1,0}(X) -g_{0,1}(Y) -\frac{1}{2} +O(C(h))$$
and $C(h) = \sum_{i=1} A_{i,1} h^i.$ By the direct computation, we can get the followings
$$\tau_{m,n} = \frac{m+n}{12mn} +O(C(h)), ~~~~~~ E_0[g_{1,0}^3 (X)] = O(C(h)),$$
$$E_0[g_{0,1}^3 (Y)] = O(C(h)), ~~~~~~ E_0[g_{1,0} (X) g_{0,1} (Y) g_{1,1} (X,Y)] = O(C(h)).$$

Then, we complete the proof of Theorem 3.2. In the similar manner of the proof of $\widetilde{M}$, we can easily obtain the Theorem 3.3.

\end{document}